\documentclass[11pt]{amsart}
\usepackage{amssymb,amscd}
\usepackage{verbatim}
\usepackage{graphicx}

\newtheorem*{rank}{Theorem}
\newtheorem*{them}{Theorem}
\newtheorem{thm}{Theorem}[section]
\newtheorem{lem}[thm]{Lemma}
\newtheorem{prop}[thm]{Proposition}

\theoremstyle{definition}
\newtheorem{defn}{Definition}[section]
\newtheorem{rem}[thm]{Remark}
\newtheorem*{remm}{Remark}

\newtheorem*{example 3.3}{Example 3.3}

\def\ov{\overline}

\newcommand{\Tors}{\operatorname{Torsion}}

\newcommand{\al}{\ensuremath{\alpha}}
\newcommand{\de}{\ensuremath{\delta}}
\newcommand{\e}{\ensuremath{\epsilon}}
\newcommand{\f}{\ensuremath{\frac}}
\newcommand{\g}{\ensuremath{\gamma}}
\newcommand{\la}{\ensuremath{\lambda}}
\newcommand{\Om}{\ensuremath{\Omega}}
\newcommand{\hra}{\hookrightarrow}
\newcommand{\x}{\ensuremath{\times}}
\newcommand{\p}{\ensuremath{\partial}}
\newcommand{\ngth}{\negthickspace}
\newcommand{\hu}[1]{\ensuremath{H_1(\underline{\ \ };#1)}}
\newcommand\lra{\longrightarrow}
\newcommand{\s}{\ensuremath{\sigma}}
\newcommand{\tha}{\ensuremath{\theta}}
\newcommand{\Th}{\ensuremath{\Theta}}

\newcommand{\2}{(2.0)-solvable}
\newcommand{\5}{(2.5)-solvable}
\newcommand{\intt}{\ensuremath{\mathrm{int}}}
\newcommand{\spin}{\ensuremath{\mathrm{spin}}}
\newcommand{\st}{\ensuremath{\mathrm{ST}}}

\newcommand{\SA}{\ensuremath{\mathcal{A}}}
\newcommand{\SC}{\ensuremath{\mathcal{C}}}
\newcommand{\SF}{\ensuremath{\mathcal{F}}}
\newcommand{\sk}{\ensuremath{\mathcal{K}}}

\newcommand{\SR}{\ensuremath{\mathcal{R}}}
\newcommand{\UU}{\mathcal{U}}
\newcommand{\NN}{\mathcal{N}}
\newcommand{\G}{\ensuremath{\Gamma}}
\newcommand{\C}{\ensuremath{\mathbb{C}}}

\newcommand{\BP}{\ensuremath{\mathbb{P}}}
\newcommand{\Q}{\ensuremath{\mathbb{Q}}}
\newcommand{\N}{\ensuremath{\mathbb{N}}}
\newcommand{\R}{\ensuremath{\mathbb{R}}}
\newcommand{\Z}{\ensuremath{\mathbb{Z}}}

\begin{document}
\baselineskip=20pt

\title{Structure in the Classical Knot Concordance Group}
\author{Tim D.\ Cochran, Kent E. Orr, Peter Teichner
\footnote{the authors were partially supported by the National
Science Foundation}}
\date{}
\maketitle

\section{Introduction}\label{sec1} Let $\SC$ be the abelian
group of topological concordance classes of knots in $S^3$.
In this paper we provide new information about its structure.
One consequence is that there is a subgroup of infinite rank
consisting entirely of knots with vanishing Casson-Gordon
invariants but whose non-triviality is detected by the
$L^{(2)}$ signature invariants of \cite{COT1}.

Recall that two knotted circles in $S^3$ are topologically
{\it concordant\/} if there is a locally flat topological
embedding of the annulus in $S^3\x[0,1]$ whose restriction to
the boundary components gives the knots. The equivalence class
of the trivial knot is the identity for $\SC$, which is a
group under connected-sum of knots, and inverses are given by
taking the mirror-image and reversing the string orientation.
A {\it slice\/} knot is one which is zero in this group, or,
equivalently bounds a locally-flat embedded disk in $B^4$.
All of our work is in this topological category.

In \cite{COT1} a geometric filtration of $\SC$ was defined
$$
0\subset\dots\subset\SF_{n.5}\subset\SF_n\subset\dots\subset
\SF_{1.5}\subset\SF_{1.0}\subset\SF_{.5}\subset\SF_0\subset\SC
$$
where $\SF_h$ consists of all {\it $h$-solvable\/} knots for $h \in
\frac{1}{2}\N_0$ (see sections~1 and 8 of \cite{COT1}). An
equivalent and more algebraic definition of these terms is reviewed in the
next section. It was shown in \cite{COT1} that $0$-solvable
knots are precisely the Arf invariant zero knots so that
$\SC/\SF_0\cong\Z_2$ given by the Arf invariant and that
$\SC/\SF_{0.5}$ is J.P.~ Levine's algebraic concordance
group which he proved was isomorphic to
$\Z^\infty\oplus(\Z_2)^\infty\oplus(\Z_4)^\infty$. It was also
shown that knots in $\SF_{1.5}$ have vanishing Casson-Gordon
invariants and thus $\SF_{1.0}/\SF_{1.5}$ has infinite rank,
as detected by Casson-Gordon invariants $[J]$ \cite{COT1}.
One of the main results of \cite{COT1} was that
$\SF_{2.0}/\SF_{2.5}$ is non-zero, in particular an infinite
set of non-slice knots with vanishing Casson-Gordon
invariants was exhibited. Here we will show that this set of
classes is linearly independent.

More precisely,

\begin{rank}\label{rank} $\SF_{2.0}/\SF_{2.5}$ has
infinite rank.
\end{rank}

The infinite set $\{K_i\}$ of non-trivial elements of
$\SF_{2.0}/\SF_{2.5}$ that was exhibited in
\cite[section~6]{COT1} was obtained from a {\it single\/}
``seedling'' ribbon knot by certain ``genetic modifications''
using a sequence of auxiliary knots $\{J_i\}$. Non-triviality
was proven by evaluating a certain higher-order invariant (an
$L^2$-signature) whose value on $K_i$ was shown to be
essentially the integral of the classical Levine-Tristram
signature function of $J_i$. This pleasing fact is verified
in Section~\ref{app} of this paper. It is simple enough to
find a set $\{J_i\}$ for which the set of these real numbers
is integrally linearly independent (and we do this in
Section~\ref{app} of the present paper). If our higher-order
invariants were additive under connected-sum, then the
Theorem would follow immediately. However, the higher-order
nature of our invariants (just as for those of Casson and
Gordon) makes it difficult even to formulate an additivity
statement. In particular, our third-order invariants depend
on choices of ``metabolizers'' (or self-annihilating
submodules) for the $1^{\text{st}}$ and $2^{\text{nd}}$-order
Blanchfield-Seifert forms, and unfortunately, the number of
such choices is usually infinite for any connected sum of
knots. To avoid the difficulties of the obvious direct
approach, we employ a slight variation which makes crucial
use of a special technical feature of the $K_i$ (arising from
the corresponding fact for the original ribbon knot), namely
that they have {\it unique\/} $1^{\text{st}}$ and
$2^{\text{nd}}$-order metabolizers.

Only the seemingly technical problem of finding a
``seedling'' ribbon knot with unique metabolizers for its
higher order linking forms (of orders $1$, $2,\dots, n-1$)
obstructs us from using the very same proof to show that
$\SF_n/\SF_{n.5}$ has infinite rank for each $n\in\Z_+$.

\begin{remm} In very recent work we have been able to show
that $\SF_{n.0}/\SF_{n.5}$ has non-zero rank for each $n\ge0$
but the proof does not seem to adapt to show infinite rank.
\end{remm}

We also include in Section~\ref{app} a proof of the following
theorem about genus one slice knots which was announced in
\cite{COT1}. It should be compared to Theorem 4 of \cite{Gi}.

\begin{them} Suppose $K$ is a 1.5-solvable knot (for example a
slice knot) whose Alexander polynomial is not 1 and which
admits a Seifert surface $F$ of genus~1. Then there exists a
homologically essential simple closed curve, $J$, on $F$ that
has self-linking number zero and such that the integral of
the Levine-Tristram signature function of $J$ vanishes.
\end{them}

We wish to thank Andrew Ranicki and Michael Larsen for helpful
contributions.

\section{$n$-solvable knots and von Neumann
$\rho$-invariants}\label{sec2}

We briefly review some of the definitions of \cite{COT1}
which are used herein.

Let $G^{(i)}$ denote the {$i$-th derived group} of a group
$G$, inductively defined by $G^{(0)}:= G$ and
$G^{(i+1)}:=[G^{(i)},G^{(i)}]$. A group is {\it $n$-solvable}
if $G^{(n+1)}=1$. For a CW-complex $W$, we define $W^{(n)}$ to
be the regular covering space corresponding to the subgroup
$\pi_1(W)^{(n)}$. If $W$ is an spin $4$-manifold then there
is an intersection form
$$
\la_n: H_2(W^{(n)})\x H_2(W^{(n)})\lra
\Z[\pi_1(W)/\pi_1(W)^{(n)}]
$$
and a self-intersection form $\mu_n$ (see \cite{Wa} chapter~5
and \cite{COT1} section~7). An {\it $n$-Lagrangian} is a
$\Z[\pi_1(W)/\pi_1(W)^{(n)}]$-submodule $L\subset
H_2(W^{(n)})$ on which $\la_n$ and $\mu_n$ vanish and which
maps (under the covering map) onto a $(1/2)$-rank direct
summand of $H_2(W;\Z)$. An {\it $n$-surface} $F$ in $W$ is a
based, immersed surface in $W$, which lifts to $W^{(n)}$. Thus
$\la_n$ and $\mu_n$ can be computed in $W$ by considering
intersections weighted by elements of
$\pi_1(W)/\pi_1(W)^{(n)}$, and an $n$-Lagrangian may be
conveniently encoded by considering a collection of
$n$-surfaces whose lifts generate it.

Suppose $K$ is a knot and $M$ is the closed $3$-manifold
resulting from $0$-framed surgery on $S^3$ along $K$.

\begin{defn}\label{n-solvable} A knot (or $M$) is {\it
$n$-solvable} ($n\in\N_0$) if $M$ bounds a spin $4$-manifold
$W$, such that the inclusion map induces an isomorphism on
first homology and such that $W$ admits two {\it dual}
$n$-Lagrangians. This means that $\la_n$ pairs the two
Lagrangians non-singularly and that their images together
freely generate $H_2(W)$. Such a $W$ is called an {\it
$n$-solution} for $K$ (or for $M$). Note that the exterior of
a slice disk is, for any $n$, an $n$-solution for the slice
knot (and for $M$) simply because the second integral homology vanishes.

A knot is {\it$(n.5)$-solvable}, $n\in\N_0$, if, in addition
to the above, one of the dual $n$-Lagrangians is the image
(under the covering map) of an $(n+1)$-Lagrangian. Then $W$
is called an {\it$(n.5)$-solution} for $K$ (or for $M$).
\end{defn}

The set $\SF_n$ (respectively $\SF_{(n.5)}$) of concordance
classes of $n$-solvable (respectively $(n.5)$-solvable) knots
is a subgroup of $\SC$.

More details can be found in \cite{COT1} sections~1, 7 and 8.

Suppose $\phi:\pi_1(M)\lra\G$ is a homomorphism. The {\it
(reduced) $L^2$-signature or von Neumann $\rho$-invariant}
$\rho(M,\phi)\in\R$ is then defined and satisfyies
$\rho(-M,\phi)=-\rho(M,\phi)$~\cite{ChG}. Recall
\cite[Definition~5.8 and Lemma~5.9]{COT1} that whenever
$(M,\phi)=\p(W,\psi)$ for some compact, oriented $4$-manifold
$W$, $\rho(M,\phi)=\s^{(2)}_\G(W,\psi)-\s_0(W)$ where
$\s^{(2)}_\G$ is the $L^2$-signature of the intersection
form on $H_2(W;\Z\G)$ and $\s_0$ is the ordinary signature.
If $\G$ is a {\it poly-torsion-free-abelian group} (henceforth
PTFA) then $\Z\G$ embeds in a skew field of fractions
$\sk_\G$ and $\s^{(2)}_\G$ may be viewed as a real-valued
homomorphism from $L^0(\sk_\G)$, so $\s^{(2)}_\G(W,\psi)$ is
a function of the Witt class of the intersection form on the
free module $H_2(W;\sk_\G)$.

We shall need only the
following properties of $\rho$ from \cite{COT1}.

\begin{enumerate}
\item[(2.2)] If $\G$ is an $n$-solvable PTFA group and $\phi$
extends over some $(n.5)$-solution $W$ for $M$, then
$\rho(M,\phi)=0$ \cite[Theorem 4.2]{COT1}. In particular, if
$K$ is a slice knot and $\phi$ extends over the exterior of a
slice disk then $\rho(M,\phi)=0$ for any PTFA group $\G$. The
reader who is not familiar with \cite{COT1} can see that the
latter follows from the very believable fact that
$H_2(W;\Z)=0$ implies that $H_2$ of the $\G$-cover of $W$ is
$\G$-torsion, hence $H_2(W;\sk_\G)=0$ \cite[Prop.~4.3]{COT1}.
\label{vanishing}
\item[(2.3)] (subgroup property) If $\phi$ factors through a
subgroup $\G'$, then $\rho(M,\phi)=\rho(M,\phi_0)$ where
$\phi_0:\pi_1(M)\lra\G'$ is the induced factorization of
$\phi$ \cite[Proposition 5.13]{COT1}. This is a consequence of
the corresponding fact for the canonical trace on a group von
Neumann algebra.\label{subgroup}
\item[(2.4)] If $\G=\Z$ and $\phi$ is non-trivial then
$\rho(M,\phi)=\int_{z\in S^1}\s(h(z))dz-\s_0(h)$ if $h$ is a
matrix representing the intersection form on
$H_2(W;\C[t,t^{-1}])/$(torsion) and $\s_0$ is the ordinary
signature \cite[Lemma 5.4, Def. 5.3]{COT1} (in fact we prove
in the Appendix that $\rho(M,\phi)$ is the integral of the
Levine-Tristram signature function of the knot $K$ normalized
to have value $0$ at $z=1$, although this more precise fact is
not strictly needed).
\label{levine}
\item[(2.5)] If $\phi$ is the trivial homomorphism then $\rho(M,\phi)=0$.
\label{map}
\end{enumerate}

If $K$ is an oriented knot then there is a {\it canonical}
epimorphism $\phi:\pi_1(M)\lra\Z$. If $K$ has Arf invariant
zero then there exists a spin $4$-manifold $W$ and a map
$\psi:\pi_1(W)\lra\Z$ such that $\p(W,\psi)=(M,\phi)$. In
fact we can assume $\pi_1(W)\cong\Z$. Such a $W$ is then a
$0$-solution for $K$ (for $M$). In particular, $\rho(M,\phi)$
is always defined for such $\phi$. Let $\rho(K)$ denote this
canonical real number (the integral of the normalized
Levine-Tristram signature function of $K$). Note
$M_{-J}=-M_J$ so $\rho(-J)=-\rho(J)$. The
proof of the following technical result is deferred to Section~\ref{app}.

\setcounter{thm}{5}
\begin{prop}\label{independence} There exists an infinite set
$\{J_i\mid i\in\Z_+\}$ of Arf invariant zero knots such that
$\{\rho(J_i)\}$ is linearly independent over the integers.
\end{prop}

\section{Constructing $n$-solvable knots by Genetic
Modification}\label{sec3} We describe a ``genetic
modification'' construction in which a given ``seedling''
knot $K$ is (sometimes subtly) altered by an auxiliary knot
$J$. We then evaluate the effect of such an alteration on the
$\rho$-invariants introduced in section~\ref{sec2}. The
construction can be summarized as follows. Seize a collection
of parallel strands of $K$ in one hand, just as you might
grab some hair in preparation for braiding. Then, treating
the collection as a single fat strand, tie it into the knot
$J$. For example, applying this to a single strand of $K$ has
the effect of altering $K$ by the addition of the local knot
$J$. This would be a rather radical alteration. Applying this
to two strands of $K$ which are part of a ``band'' of a
Seifert surface for $K$ has the effect of tying that band
into a knot $J$. For this simple case the construction
agrees with the one used often by A. Casson, C. Gordon, P.
Gilmer, D. Ruberman, C. Livingston and others to create
knots with identical Alexander modules but different
Casson-Gordon invariants. However, in the more complicated
cases, the genetic modification construction can be seen to
be somewhat ``orthogonal'' to that of the construction of
Casson-Gordon and to the more recent modification
constructions of T. Stanford and K. Habiro. Our construction,
which had its origins in \cite{CO}, is more similar to J.
Levine's realization constructions via surgery on links.
Various versions of this construction have been useful in
knot and link theory (\cite{CO}, \cite[section~6]{COT1}).
More details and variations on the construction appear in
\cite{C}. In the application relevant to this paper, we will
choose $K$ to be a ribbon knot and choose a circle in
$S^3\backslash K$ which lies deep in the derived series of
$\pi_1(S^3\backslash K)$ and yet bounds an embedded disk in
$S^3$. The knot will pierce this disk many times. The
alteration is to cut open $K$ along this disk and tie all the
strands passing through into a knot $J$ as shown in the
right-most part of Figure 1. Details of the general
construction follow.

\begin{figure}[ht]
\begin{center}
        \includegraphics[scale=.6]{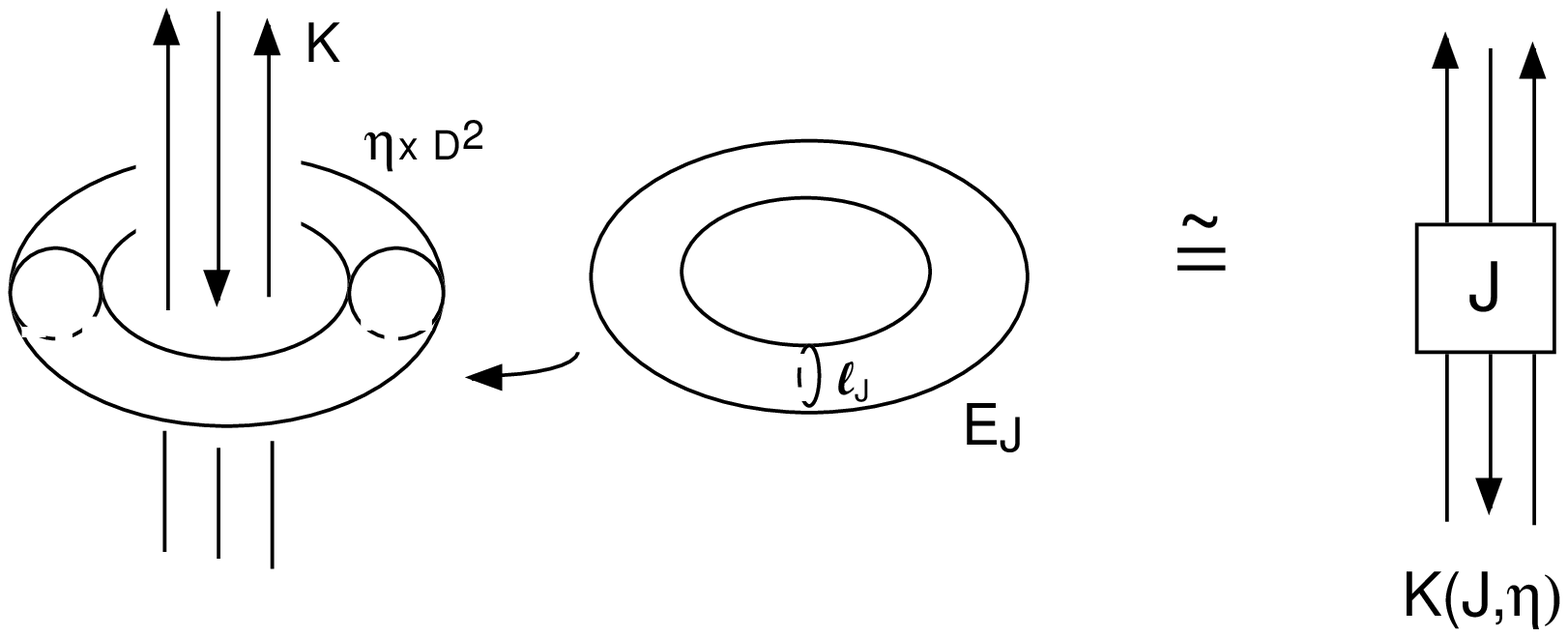}
\caption{}\label{fig1}
\end{center}
\end{figure}

Let $M$ and $M_J$, respectively, denote the zero framed
surgeries on the knots $K$ and $J$, and $E$ and $E_J$ denote
their exteriors. Suppose $\eta$ is an oriented simple closed
curve in $E$, which is unknotted in $S^3$. Choose an
identification of a tubular neighborhood of $\eta$ with
$\eta\x D^2$ in such a way that $\eta\x\{1\}\subset\eta\x\p
D^2$ is a longitude $\ell_\eta$, and $\{*\}\x\p D^2$ is a
meridian $\mu_\eta$. Form a new oriented manifold
$E'=(E-\intt(\eta\x D^2))\cup E_J$ by an identification of
$\p E_J$ with $\eta\x\p D^2$ which sends
$\mu^{-1}_\eta$ to the longitude of $J$, denoted $\ell_J$, and
sends $\ell_\eta$ to $\mu_J$. Note that $\p E'=\p E=K\x D^2$
and since $E\cup(K\x D^2)=S^3$,
$E'\cup(K\x D^2)=(S^3-\intt(\eta\x D^2))\cup E_J$. Since
$\eta$ is unknotted, $S^3-\intt(\eta\x D^2)$ is a solid torus
$\st$ and $\st\cup E_J\cong S^3$ as can be confirmed by
checking the identifications. Hence $E'\cup(K\x D^2)\cong
S^3$. Therefore $E'$ is the exterior $S^3\setminus K'$ of a
knot $K'$ which is the image of $K$ under the identification
$E'\cup(K\x D^2)\cong S^3$. The new knot $K'=K(J,\eta)$ is a
function of $K$, $J$ and $\eta$ as well. Since $K\subset\st$,
$K'$ is a satellite of $J$. It is left to the reader to see
that $K'$ is indeed the result of tying the strands of $K$
that ``pass through the $2$-disk spanned by $\eta$'' into the
knot $J$. Finally, if we let $M'$ denote the zero framed
surgery on $K'$, then $M'=(M-\intt(\eta\x D^2))\cup E_J$. We
shall show the following:

\begin{prop}\label{construction} If $K$ is $n$-solvable,
$\eta\in\pi_1(M)^{(n)}$ and $J$ has Arf invariant zero, then
$K'=K(J,\eta)$ is $n$-solvable.
\end{prop}

Hence one speculates that if one begins with a slice knot $K$,
which is of course $n$-solvable for each $n$, and
$\eta\in\pi_1(M)^{(n)}-\pi_1(M)^{(n+1)}$, then many different
$n$-solvable knots could be constructed as long as the derived
series of $\pi_1(M)$ does not stabilize. In fact it is proved
in \cite{C} that the derived series of a knot group cannot
stabilize unless it has Alexander polynomial $1$ and that of
$\pi_1(M)$ cannot stabilize unless the Alexander polynomial
has degree~$0$ or $2$. The above speculation is also
confirmed in \cite{C}.

We remark that if $n$-solvable is replaced by rationally
$n$-solvable (see \cite[section~4]{COT1}) then the analogous
result holds without a condition on the Arf invariant of $J$.
This is of interest in studying knots which bound disks in
rational homology balls.

We now want to show that many of these families of knots are
different, even up to concordance. For this purpose we
consider the $L^2$-signature invariants of section~\ref{sec2}
\cite[section~5]{COT1}.

Suppose one is given homomorphisms
$\phi:\pi_1(M)\to\G$, $\phi_J:\pi_1(M_J)\to\G$ such that
$\phi([\eta])=\phi_J([\mu_J])$ and $\G$ is a PTFA group.
Then, if $M'$ is as above, a unique $\phi':\pi_1(M')\to\G$ is
induced which extends $\phi$ and $\phi_J$.

\begin{prop}\label{additivity} Given $K$, $J$, $\eta$, $\phi$,
$\phi_J$ as above,\newline
$\rho(M',\phi')=\rho(M,\phi)+\rho(M_J,\phi_J)$
whenever the right hand side is defined.
\end{prop}

\begin{proof}[Proof of Proposition \ref{construction}] Let $W$
be an $n$-solution for $K$ and let $W_J$ be the $0$-solution
for $J$ with $\pi_1\cong\Z$ discussed above
Proposition~\ref{independence}. So $\p
W_J=M_J=E_J\cup(S^1\ngth\x\ngth D^2)$ where
$S^1\ngth\x\ngth\{*\}$ is $\mu_J$ and $\{*\}\ngth\x\ngth\p
D^2$ is $\ell_J$. Let $W'$ be the $4$-manifold obtained from
$W_J$ and $W$ by identifying the solid torus $S^1\ngth\x\ngth
D^2\subset\p W_J$ with $\eta\ngth\x\ngth D^2\subset\p W$.
Observe that $\p W'=M'$, the zero surgery on $K'=K(J,\eta)$.
We claim that $W'$ is an $n$-solution for $K'$. First consider
the Mayer-Vietoris sequence below
$$
0\lra H_2(W)\oplus H_2(W_J)
\xrightarrow{\pi_*} H_2(W')\xrightarrow{\p_*} H_1(W\cap W_J)
\xrightarrow{i_*}H_1(W)\oplus H_1(W_J)
$$
Since $W\cap W_J\simeq S^1$, $H_2(W\cap W_J)=0$ for any
coefficients. Since the inclusion $W\cap W_J\to W_J$ induces
an isomorphism on $\pi_1$, the map $i_*$ is a monomorphism
with any coefficients. Thus $\pi_*$ is an isomorphism with any
coefficients, and the intersection and self-intersection
forms on $H_2$ split naturally. We may think of an
$n$-Lagrangian with its $n$-duals for $W$ as being generated
by finite collections of based surfaces in $W$ each of which
lifts to the $\pi_1(W)^{(n)}$-cover (these were called
``$n$-surfaces'' in \cite{COT1}; sections~7--8). These same
surfaces are clearly $n$-surfaces in $W'$ since
$n^{\text{th}}$-order commutators in $\pi_1(W)$ are
$n^{\text{th}}$-order commutators in $\pi_1(W')$. Similarly
consider the collection of ``$0$-surfaces'' generating a
$0$-Lagrangian and its duals for $W_J$. Since the map
$\pi_1(W_J)\to\pi_1(W')/\pi_1(W')^{(n)}$ is trivial (since
$\pi_1(W_J)$ is generated by $\eta$), these $0$-surfaces are
$n$-surfaces in $W'$. It then follows easily by naturality
that the union of these collections of
$n$-surfaces constitutes an $n$-Lagrangian with $n$-duals for
$W'$ (see \cite{COT1}; sections~7--8).
\end{proof}

\begin{proof}[Proof of Proposition \ref{additivity}]
Suppose $\p(W,\psi)=(M,\phi)$ and
$\p(W_J,\psi_J)=(M_J,\phi_J)$ for some $H_1$-bordisms as
described in Section~\ref{sec2}. Note that $M_J=E_J\cup(S^1\x
D^2)$ where $\{1\}\x\p D^2$ is a longitude of $J$. Let $W'$
be the $4$-manifold obtained from $W_J$ and $W$ by identifying
the solid torus $S^1\ngth\x\ngth D^2\subset\p W_J$ with
$\eta\ngth\x\ngth D^2\subset\p W$. Observe that $\p W'=M'$ and
that $\psi$ and $\psi_J$ piece together to give an extension
of $\phi'$ to $\psi':\pi_1(W')\to\G$. Thus $(W',\psi')$ can be
used to compute $\rho(M',\phi')$. It now suffices to show that
the natural inclusions induce an isomorphism $H_2(W)\oplus
H_2(W_J)\to H_2(W')$ with $\sk$ coefficients. As above,
$H_2(W\cap W_J)=0$ with any coefficients. Now, appealing to
the Mayer-Vietoris sequence, it will suffice to show that
$H_1(W\cap W_J;\sk)\to H_1(W_J;\sk)$ is a monomorphism. If
$\phi([\mu_J])=\phi([\eta])\neq0$ then
$H_1(W\cap W_J;\sk)=0$ since the induced $\G$-cover of a
circle is a union of lines (alternatively use \cite{COT1};
Proposition~2.11). If $\phi([\mu_J])=0$ then $\phi_J$ is the
zero map since $\pi_1(M_J)$ is normally generated by
$[\mu_J]$. Hence we may assume $\psi_J$ is the zero
homomorphism. Therefore the coefficient system
$\psi_J:\Q[\pi_1(W_J)]\to\sk$ factors as the augmentation
$\e:\Q[\pi_1(W_J)]\to\Q$ followed by the inclusion
$\Z=\Q[\{1\}]\hra\Q\G\hra\sk$. Thus
$H_1(W_J;\sk)\cong H_1(W_J;\Q)\otimes_\Q\sk\cong\sk$ and
similarly for $H_1(W\cap W_J;\sk)$. Since the inclusion
$W\cap W_J\to W_J$ induces an isomorphism on \hu{\Q}, it
induces an isomorphism on \hu{\sk} as well.
\end{proof}

The following application is what we will use {\it in the proof of}
the main theorem.

\begin{example 3.3} Suppose $K$ is a ribbon knot,
$\eta\in\pi_1(M)^{(n)}$ and
$J$ has Arf invariant zero. Then $K'=K(J,\eta)$ is
$n$-solvable by Proposition~\ref{construction}. Let $W$ be
the exterior of a ribbon disk for $K$ and let $W_J$ be the
$(0)$-solution for $J$ with $\pi_1(W_J)\cong\Z$ as in
section~\ref{sec2}. Then let $W'$ be the $n$-solution for
$K'$ formed as in the proof of Proposition~\ref{construction}
by gluing $W_J$ to $W$ along $\eta\x D^2$. Suppose
$\psi':\pi_1(W')\to\G$ is a homomorphism defining, by
restriction, $\phi'$, $\psi$, $\phi_J$ and $\psi_J$ from
(respectively) $\pi_1(M')$, $\pi_1(W)$,
$\pi_1(M_J)$, and $\pi_1(W_J)$. Then
$\rho(M',\phi')=\rho(M_J,\phi_J)$ by
Proposition~\ref{additivity} and 2.2. Then, (since
$\pi_1(W_J)\cong\Z$ is generated by $\eta$) using 2.3, 2.4,
and 2.5, $\rho(M',\phi')$ equals $\rho(J)$ if
$\phi(\eta)\neq1$ and equals $0$ if $\phi(\eta)=1$.
\end{example 3.3}

\section{The Main Theorem}

\begin{thm}\label{main theorem} $\SF_{(2.0)}/\SF_{(2.5)}$ has
infinite rank.
\end{thm}

\begin{proof}  It is sufficient to exhibit an infinite set of
\2 knots $K_i$, $i\in\Z_+$, such that no non-trivial linear
combination is \5. Let $K_r$ be the ``seedling'' ribbon knot
shown in Figure~3.1 and $\eta$ be the designated circle just
as was used in \cite[section~6]{COT1}. Let $J_i$ be the knots
of Proposition~\ref{independence} and let $K_i=K_r(J_i,\eta)$
be the family of knots resulting from the grafting
construction of Example~3.3. It was shown in
\cite[section~6]{COT1} that each of these knots is 2-solvable
and not \5. Suppose that a non-trivial linear combination
$\#^m_{i=1}n'_iK_i$, $n'_i\neq0$, were \5. We shall derive a
contradiction.

We may assume all $n'_i>0$ by replacing $K_i$ by $-K_i$ if
$n'_i<0$. We may also assume that if $m=1$ then $n'_1>1$. For
simplicity let $M_i$ denote $M_{K_i}$, and note that
$-M_i=M_{-K_i}$ because $-K_i$ can be obtained by applying a
reflection to $(S^3,K_i)$. Let $M_0$ denote the zero surgery
on $\#n'_iK_i$ and let $W_0$ denote the putative
(2.5)-solution.

\begin{figure}[ht]
\begin{center}
        \includegraphics[scale=.6]{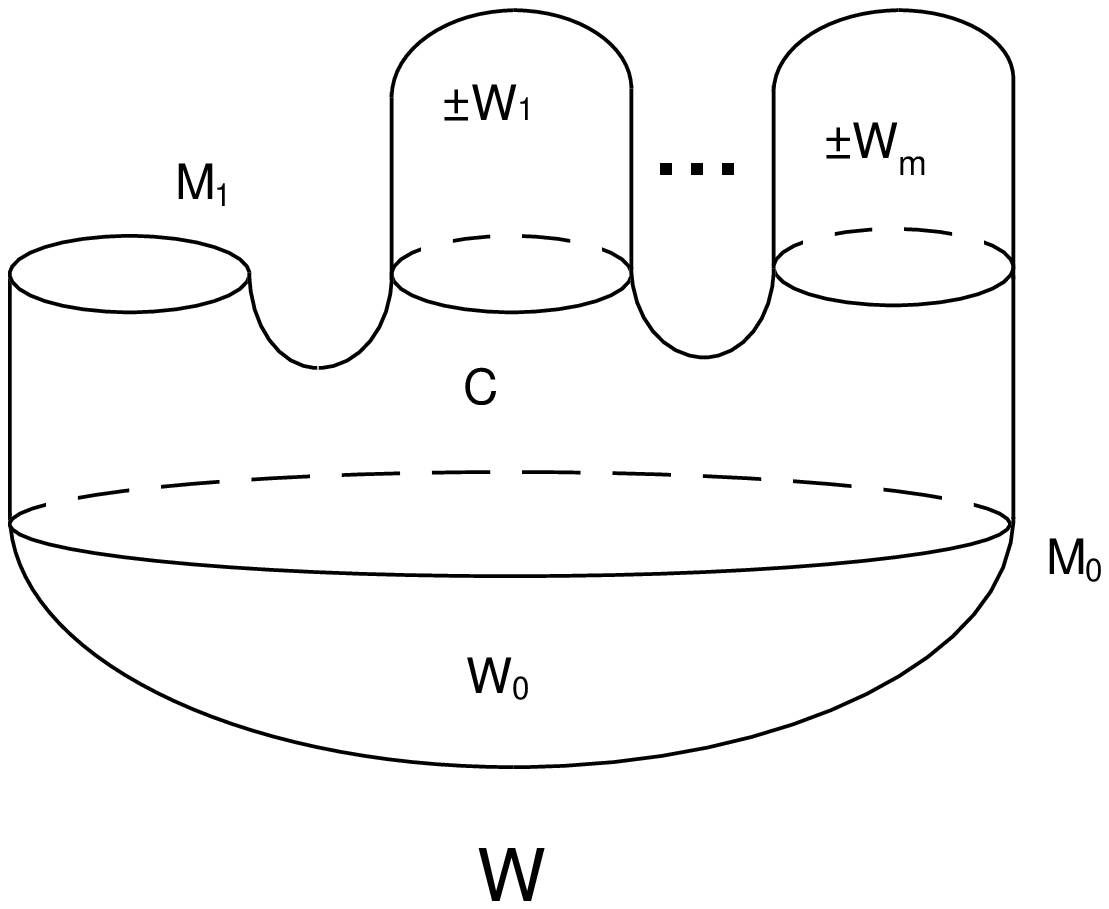}
\caption{}
\end{center}
        \end{figure}

\begin{figure}[ht]
\begin{center}
       \includegraphics[scale=.7]{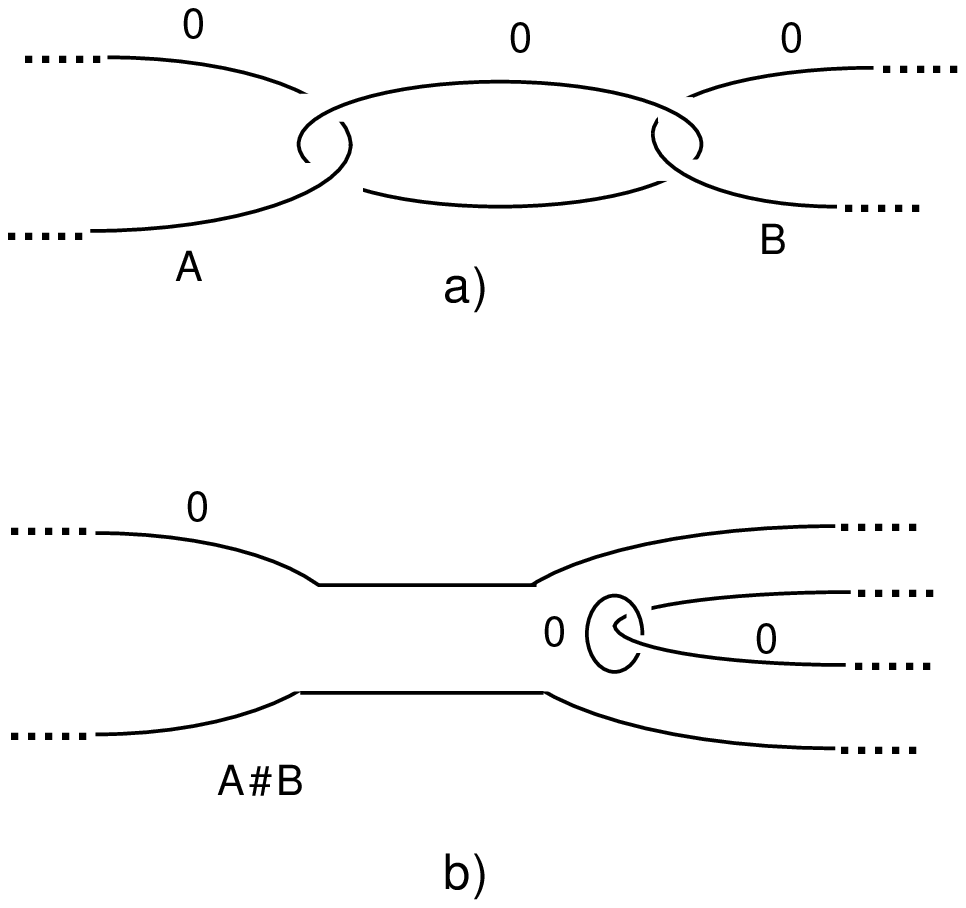}
\caption{}
\end{center}
        \end{figure}

\medskip
  From these assumptions we construct a very specific
2-solution $W$ for $M_1$. This is shown schematically in
Figure~2 where $C$ and $W_i$ are to be defined. The $W_i$,
$i>0$, are the specific $2$-solutions for $M_i$ constructed
as in Example~3.3.  There are $n_i$ copies of $W_i$ where
$n_1=n'_1-1$ and $n_i=n'_i$ if $i>1$. The $4$-manifold $C$ is
merely a standard cobordism between zero surgery on a
connected sum of knots and the disjoint union of the zero
surgeries on its summands. For the case of just two knots
$A\#B$, the manifold $C$ is described as follows. Beginning
with a collar on $M_A\amalg M_B$, add a 1-handle to get a
connected $4$-manifold whose ``upper'' boundary is given by
surgery on the split link $A\amalg B\hookrightarrow S^3$, each
with zero framing. Next add a zero framed 2-handle along a
circle as shown in Figure~3a. This completes the
description of $C$ in this simple case. We need only show that
the $3$-manifold in Figure~3a is homeomorphic to $M_{A\#B}$.
This is accomplished most easily in the language of Kirby's
calculus of framed links \cite{Ki}. By ``sliding'' the $A$
circle over the $B$ circle, one arrives at the equivalent
description shown in Figure~3b. But now the circle labeled
$B$ may be canceled by the small linking circle, leaving
only the desired zero surgery on the connected sum. By
iterating this idea, one sees that our $C$ has a handlebody
decomposition, relative to $\coprod^m_{i=1}n'_iM_i$,
consisting of $(\sum^m_{i=1}|n_i|)$ 1-handles and the same
number of 2-handles. The 1-handles have no effect on $H_1$ or
$H_2$, while the 2-handles serve to equate all the meridional
generators of $H_1$ (and thus do not affect $H_2$). Hence
$H_1(C;\Z)\cong\Z$ and the inclusion from any of its boundary
components induces an isomorphism on $H_1$; and $H_2(C)\cong
H_2(\amalg n'_iM_i)$. Moreover, in $H_2(C)$, the generator of
$H_2(M_0)$ is equal to the sum of generators of $H_2$ of the
other components of $\p C$.

We now verify that $W$ is an $H_1$-bordism for $M_1$, and
identify $H_2(W)$. The inclusions induce isomorphisms
$H_1(M_i;\Z)\to H_1(W_i;\Z)$. It then follows that
$H_1(M_1)\to H_1(W)\cong\Z$ is an isomorphism. Now consider a
Mayer-Vietoris sequence for $W=C\cup W^*$ where
$W^*=\ov{C-W}$. We see that $H_1(C\cap W^*)\to H_1(W^*)$ is
injective by the remarks above. Note that the boundary map
$H_3(W_i,M_i)\to H_2(M_i)$ is an isomorphism since it is dual
to the inclusion $H^1(W_i)\lra H^1(M_i)$. Thus for $i\ge0$,
$H_2(M_i)\to H_2(W_i)$ is the zero map. Therefore $H_2(C\cap
W^*)\to H_2(W^*)$ is the zero map. By the last sentence of
the previous paragraph, $H_2(C\cap W^*)\to H_2(C)$ is
surjective. It follows that
$H_2(W)\cong H_2(W^*)\cong H_2(W_0)\bigoplus^m_{i=1}n_i
H_2(W_i)$. It is not difficult to see that $W$ is a {\it
spin} bordism since each individual piece is spin with 2 spin
structures and $\Om^{\spin}_3(\Z)\cong\Z_2$ is given by the
Arf invariant.

To show that $W$ is a 2-solution for $M_1$, we must exhibit a
2-Lagrangian with 2-duals. But this is obtained merely by
taking the ``union'' of the 2-Lagrangians and 2-duals for
$W_0$ and each $\pm W_i$ which appears as part of $W$.
(Recall that since $W_0$ is a (2.5)-solution it is also a
2-solution). More precisely suppose, for example, that
$\{L_1,\dots,L_m\}$, $\{D_1,\dots,D_m\}$ are 2-surfaces in
$W_0$ which generate the 2-Lagrangian and its dual
2-Lagrangian for $W_0$. In particular these surfaces lift to
$W^{(2)}_0$ and so the image of $\pi_1(L_i)$ in $\pi_1(W_0)$
is contained in $\pi_1(W_0)^{(2)}$. Thus this image is
contained in $\pi_1(W)^{(2)}$ and so these surfaces lift to
$W^{(2)}$. Similarly by functoriality of intersection with
twisted coefficients, these surfaces have the required
intersection properties when considering the intersection
form on $W$ with $\Z[\pi_1(W)/\pi_1(W)^{(2)}]$ coefficients.
An identical argument is used for each $W_i$, making it clear
that the union of these 2-surfaces represents a 2-Lagrangian
and 2-duals for $W$, completing the demonstration that $W$ is
a 2-solution for $M_1$.

Since $W$ is a 2-solution for $M_1$, Theorem 4.6 of
\cite{COT1} guarantees the existence of certain non-trivial
homomorphisms $\phi_2:\pi_1(M_1)\to\G^U_2$, for a certain {\it
universal solvable group} $\G^U_2$, which extend to
$\pi_1(W)$. Moreover, if $W$ is fixed, such homomorphisms
actually factor through a much smaller group (the image of
$\pi_1(W)$, for example). This improvement was mentioned in
Remark~4.7 of \cite{COT1} and was discussed in detail in
section~6 of that paper for precisely the case at hand. We
shall repeat some of that argument. Let $\G_0=\Z$ and let
$\phi_0:\pi_1(M_1)\to\G_0$ be the canonical epimorphism which
extends uniquely to an epimorphism $\psi_0:\pi_1(W)\to\G_0$.
Recall that the classical Alexander module
$\SA_0(M_1)=H_1(M_1;\Q[t,t^{-1}])$ is isomorphic to
$\Q[t,t^{-1}]/(p(t))^2$ where $p(t)=t^{-1}-3+t$ (this
computation was discussed, but left to the reader in section~6
of \cite{COT1}). Let $\SA_0(W)=H_1(W;\Q[t,t^{-1}])$. By
\cite[Theorem~4.4]{COT1}, since $W$ is a 1-solution for $M_1$,
the kernel of the inclusion-induced map
$j_*:\SA_0(M_1)\to\SA_0(W)$ is self-annihilating with respect
to the Blanchfield form $B\ell_0$. Since $\SA_0(M_1)$ has a
unique proper submodule $P_0$, the latter is in fact this
kernel. Choose a non-zero $p_0\in P_0$, inducing
$\phi_1:\pi_1(M_1)\to\G^U_1$ by \cite[Theorem~3.5]{COT1}
(recall $\G^U_1=\Q(t)/\Q[t,t^{-1}]\rtimes\G_0$). By
\cite[Theorem~3.6]{COT1} $\phi_1$ extends to
$\psi_1:\pi_1(W)\to\G^U_1$. Using the argument of
(\cite{COT1}, see just before Proposition~6.1), we can
replace $\G^U_1$ by a subgroup $\G_1$ which contains the
image of $\phi_1$ and is isomorphic to
$\Q[t,t^{-1}]/(p(t)^m)\rtimes\G_0$ for some positive integer
$m$; replace $\phi_1$ by restricting its image and replace
$\psi_1$ by a new map extending this restriction. We
re-label these new maps by $\phi_1$ and $\psi_1$. Continuing
as in \cite{COT1}, we choose a subring $\SR_1$ of the field
of fractions of $\Z\G_1$ where
$\SR_1=(\Q[[\G_1,\G_1]]-\{0\})^{-1}\Q\G_1$ and then set
$\SA_1(M_1)=H_1(M_1;\SR_1)$, $\SA_1(W)=H_1(W;\SR_1)$ using
the coefficient systems $\phi_1$ and $\psi_1$ respectively.
Then, since $W$ is a 2-solution, the kernel of
$j_*:\SA_1(M_1)\to\SA_1(W)$ is self-annihilating with respect
to the non-singular linking form $B\ell_1$
\cite[Theorem~4.4]{COT1}. By Proposition~6.1 of \cite{COT1},
$\SA_1(M_1)$ has a unique self-annihilating submodule $P_1$
which is therefore the kernel of $j_*$. Choose a non-zero
element $p_1\in P_1$ (by \cite[Proposition~6.2b]{COT1} $P_1$
is non-trivial). This induces $\phi_2:\pi_1(M_1)\to\G_2$
($\G_2=\sk_1/\SR_1\rtimes\G_1$ where $\sk_1$ is the quotient
field of $\SR_1$). We note that the loop $\eta$ is chosen so
that $\phi_2(\eta)\neq e$ (see below Proposition~6.1 in
\cite{COT1}). Since $W$ is a $2$-solution for $M_1$,
Theorem~3.6.1 of \cite{COT1} applies (with $n=2$, $x=p_1$) to
show that $\phi_2$ extends to $\psi_2:\pi_1(W)\to\G_2$.
Therefore $\rho(M_1,\phi_2)$ is defined and can be computed
using $(W,\psi_2)$.

We shall now compute $\rho(M_1,\phi_2)$ using $(W,\psi_2)$.
Let $\phi_{(i,j)}$ denote the restriction of $\psi$ to the
$j^{\text{th}}$ copy of $\pi_1(M_i)$ $1\le j\le n_i$, and let
$\phi_0$ denote the restriction of $\psi$ to $\pi_1(M_0)$.
Note that each of these homomorphisms is non-trivial since the
generator of $H_1(C;\Z)\cong H_1(M_1;\Z)\cong\Z$ is carried
by each $M_i$ and since $\phi_2$ and $\psi_2$ agree on the
abelianizations. Let $\sk_2$ denote the quotient (skew) field
of $\Z\G_2$. Consequently $H_*(M_i;\sk_2)=0$ for $i\ge0$
(Propositions~2.9 and 2.11 of \cite{COT1}), and hence a
Mayer-Vietoris sequence shows that
$H_2(W;\sk_2)\cong H_2(W_0;\sk_2)\oplus H_2(C;\sk_2)\oplus
H_2(W_1;\sk_2)\dots\oplus H_2(W_m;\sk_2)$ where $-W_i$ occurs
$n_i$ times and the coefficient systems on the subspaces of
$W$ are induced by inclusion. Similarly the intersection form
on $H_2(W;\sk_2)$ splits as such a direct sum. Moreover we
claim that $H_2(C;\sk_2)=0$. Let $\p_+ C=\p W_0$ and $\p_-
C=\coprod n'_iM_i$. We have observed that $(C,\p_-C)$ is a
relative 2-complex with an equal number of 1-handles and
2-handles. The claim will follow from Lemma~\ref{bordism}
which shows that, even though $C$ is not an integral homology
cobordism, it is a $\sk_2$-homology cobordism.

\begin{lem}\label{bordism} Suppose $(C,\p C)$ is a compact,
oriented $4$-dimensional Poincar\'e complex such that $\p
C=\p_+C\amalg\p_-C$, $(C,\p_-C)$ is homotopy equivalent to a
finite (relative) $2$-complex with no $0$-handles and a equal
number of $1$- and $2$-handles. Suppose also that
$\beta_1(\p_+C)=1$ and that $\phi:\pi_1(C)\to\G$ is
non-trivial on $\pi_1(\p_-C)$ where $\G$ is a
poly-torsion-free-abelian group with quotient field $\sk_\G$.
Then
$$
H_*(C,\p_+C;\sk_\G)\cong H_*(C,\p_-C;\sk_\G)\cong
H_*(C;\sk_\G)\cong0.
$$
\end{lem}

\begin{proof} It follows that $\phi$ is non-trivial on
$\pi_1(C)$ and $\pi_1(\p_+C)$ and so
$H_0(C;\sk)\cong H_0(\p_+C;\sk)\cong H_0(\p_-C;\sk)\cong
H_0(C,\p_+C;\sk)\cong H_0(C,\p_-C;\sk)\cong0$ by
\cite[Proposition~2.9]{COT1}. Moreover $H_1(C,\p_+C;\sk)\cong
H^3(C,\p_-C;\sk)\cong0$ since $(C,\p_-C)$ is a $2$-complex.
Since $\beta_1(\p_+C)=1$, $H_1(\p_+C;\sk)=0$ by
\cite[Proposition~2.11]{COT1} and $H_2(\p_+C;\sk)\cong
H^1(\p_+C;\sk)\cong0$ by Remark~2.8.3 of that paper.
Combining this with the previous facts, we see that
$H_1(C;\sk)\cong0$ and hence that $H_1(C,\p_-C;\sk)\cong0$ by
consider the sequence of the pair $(C,\p_-C)$. Finally, note
that the chain complex obtained from the cell-structure  of
the $2$-complex $(C,\p_-C)$, by lifting to the $\G$-cover and
tensoring with $\sk$, has only two non-zero terms, which are
free $\sk$-modules of the same rank. Since $H_1$ of this
chain complex is zero, the boundary map $\p_2$ is an
epimorphism and hence an isomorphism. Thus
$H_2(C,\p_-C;\sk)=0$ which implies $H_2(C;\sk)$ vanishes,
implying $H_2(C,\p_+C;\sk)=0$, and the claimed results in
relative homology and hence homology now follow.
\end{proof}

Thus, since $\s^{(2)}$ is a homomorphism on Witt classes of
non-singular forms with $\sk_2$ coefficients, we see that
$$
\rho(M_1,\phi_2) = \rho(M_0,\phi_0) + \sum^m_{i=1}
\sum^{n_i}_{j=1}\rho(-M_i,\phi_{(i,j)}).
$$
Since $\phi_0$ extends to the (2.5)-solution $W_0$,
$\rho(M_0,\phi_0)$ vanishes by (2.2). Let $\e_{ij}$ equal
$0$ or $1$ according as $\phi_{(i,j)}(\eta)$ equals $0$ or
not. Since $\phi_{(i,j)}$ extends to $-W_i$, Example~3.3
establishes that
$\rho(-M_i,\phi_{(i,j)})=-\e_{ij}\rho(J_i)$. This shows that
$\rho(M_1,\phi_2)+\sum^m_{i=1}c_i\rho(J_i)=0$ where
$c_i=\sum^{n_i}_{j=1}\e_{ij}$ is non-negative.

It will now suffice to show that $\rho(M_1,\phi_2)$ equals
$\rho(J_1)$, for this will, for any $c_i$, contradict
Proposition~\ref{independence}. The argument of
\cite[section~6]{COT1} (outlined earlier in this proof as
regards extension to $W$) shows that $\phi_2$ extends to the
manifold $W_1$ as constructed in Example~3.3, the crucial
facts being that $W_1$ is a $2$-solution and the uniqueness of
the self-annihilating submodules for the ordinary and
first-order generalized Alexander modules of $K_1$. Hence, by
2.3 $\rho(M_1,\phi_2)=\rho(J_1)$ since $\phi_2(\eta)\neq
e$. This contradiction establishes that no non-trivial linear
combination of the $K_i$ is \5 and hence proves the Theorem.
\end{proof}

\section{Appendix: Some Results on Slice Knots}\label{app}

The purpose of this section is to prove the following two
results, as well as to prove Proposition~\ref{independence}.
Only the latter is required for our main theorem.

\begin{prop} \label{LT}
For $\omega\in S^1$, let $\sigma_\omega$ be the
Levine-Tristram signature of a knot $K$. Then the (reduced)
$L^2$-signature of $K$ (denoted $\rho(K)$ in
Section~\ref{sec2}) is the integral of these signatures
$\sigma_\omega$, integrated over the circle normalized to
length one.
\end{prop}

\begin{thm}\label{genus one} Suppose $K$ is a 1.5-solvable
knot (for example a slice knot) whose Alexander polynomial is
not 1 and which admits a Seifert surface $F$ of genus~1. Then
there exists a homologically essential simple closed curve,
$J$, on $F$ that has self-linking number zero and such that
the integral of the Levine-Tristram signature function of $J$
vanishes.
\end{thm}

The above theorem, announced in \cite{COT1} should be
compared to \cite[Theorem 4]{Gi}. It is not known if that
theorem implies ours or vice-versa!

\begin{proof}[Proof of Theorem \ref{genus one}] Suppose
$W$ is a 1.5-solution for $M$, the zero framed surgery on
$K$. Let $G=\pi_1(W)$, let
$\frak A=G^{(1)}/G^{(2)}\otimes_{\Z[t,t^{-1}]}\Q[t,t^{-1}]$, and
let $\G= \frak A\rtimes\Z$. Since $\frak A$ is a torsion-free abelian
group, $\G$ is a PTFA (poly-torsion-free-abelian-group) and
is 1-solvable. Moreover there is a canonical homomorphism
$\psi:G\to\G$ since
$G/G^{(2)}\cong(G^{(1)}/G^{(2)})\rtimes\Z$. Let
$\phi:\pi_1(M)\to\G$ be the composition of
$j_*:\pi_1(M)\to\pi_1(W)$ with $\psi$. By 2.2,
$\rho(M,\phi)=0$. Since $H_1(M)\cong H_1(W)\cong
H_1(\G)\cong\Z$ we can consider the inclusion induced map on
infinite cyclic covers,
$j_*:H_1(M;\Q[t,t^{-1}])\to H_1(W;\Q[t,t^{-1}])$. The former
group is isomorphic to the classical rational Alexander module
$\SA$ of $K$ (since the longitude lies in the second derived)
and the latter is the group denoted $\frak A$ above. Since $W$ is a
1-solution for $M$, the kernel of $j_*$ is a submodule
which is self-annihilating with respect to the classical
Blanchfield form (apply \cite[Theorem 4.4]{COT1} with $n=1$
and $\G=\Z$). This fact is well known in case $W$ is the
exterior of a slice disk for $K$.  In particular, this implies that the
kernel of $j_*$ and the image of $j_*$ have rank (over $\Q$) equal to
one-half the degree of the Alexander polynomial, which is, by assumption,
positive. Moreover this implies
that $K$ is algebraically slice. The following algebraic fact
about a genus 1 algebraically slice knot can be found in
\cite[Section 5]{Gi}. $\SA$ is a cyclic module with precisely 2 proper
submodules each generated by a simple closed curve on $F$ of zero
self-linking number. Let
$J$ denote the proper submodule contained in kernel $(j_*)$. It follows
that $J$ represents an element of kernel $\phi$. Now
construct a cobordism over $\G$ from $(M,\phi)$ to
$(M_J,\phi')$ as follows, where $M_J$ is zero surery on $J$. To
$M\x\{1\}\subset M\x[0,1]$ attach a zero-framed 2-handle along
$J\x\{1\}$. The result of such a surgery can be seen to be homeomorphic to
$M_J\#S^1\x S^2$ by sliding the zero-framed 2-handle
$K\x\{1\}$ twice over $J$ so that it becomes unknotted and
unlinked from $J$. Adding a 3-handle along $\{pt\}\x S^2$
completes the cobordism. Since $J\subset\ker\phi$, $\phi$
extends over the cobordism, inducing $\phi':\pi_1(M_J)\to\G$.
Since $\pi_1(M_J)$ is normally generated by the meridian of
$J$, which clearly is null-homologous in the cobordism,
$\phi'$ factors through $[\G,\G]=\frak A$. Since $H_1(M_J)\cong\Z$,
the image of $\phi'$ is either $0$ or $\Z$.  We claim it is the latter.
For since the image of $j_*$ is not zero, there is an element of the
Alexander module that maps non-trivially under $j_*$. Let $\gamma
\in \pi_1(M-F)$  represent this element. Then $\gamma$ maps
non-trivially under $\phi$. Since the $2$ and $3$-handles above are
attached disjointly from $M-F$, $\gamma$ represents an element of
$\pi_1(M_J)$ which maps non-trivially under $\phi'$. Thus the image of
$\phi'$ is $Z$.  By Lemma~\ref{bordism}, this cobordism is a
$\sk_\G$-homology cobordism and thus
$0=\rho(M,\phi)=\rho(M_J,\phi')$. Moreover
$\rho(M_J,\phi')=\rho(J)$ which equals the integral of the
signature function of $J$ (by 2.3, 2.4 and
Proposition~\ref{LT}).
\end{proof}

We now move to the proof of Proposition~\ref{LT}. Note that
even though the $\sigma_\omega$ are integers, the integral is
in general only a real number. We first give the relevant
definitions and then derive Proposition~\ref{LT} from
Lemma~\ref{lem:F=W}.

First recall that the Levine-Tristram signature is defined as
follows: Pick a Seifert surface $F$ for the knot, together
with a basis of embedded curves $a_1,\dots,a_{2g}$ on
$F$. Using the positive push-off's $a_i^\uparrow$ of $a_i$
into 3-space, the corresponding Seifert matrix is defined using linking
numbers in $S^3$:
$$
S_{ij} := lk(a_i, a_j^\uparrow).
$$
Then for $ \omega \in S^1,  \sigma_\omega$ is the signature of the complex
hermitian form
$$
\lambda_\omega(F):=(1- \omega^{-1} ) \cdot S + (1-
\omega)\cdot S^T.
$$

Instead of choosing a Seifert surface for the knot $K$, one
can also choose a $4$-manifold $W$  which bounds the
$0$-surgery $M_K$ on $K$. One can also arrange that $W$ has
signature zero (by adding copies of $\pm\C\BP^2$) and that
$ \pi_1W \cong\Z$, generated by a meridian $t$ of the knot.
This follows from the fact that the bordism group
$\Omega_3(S^1)$ vanishes.

For the purpose of this section, let's call a $4$-manifold $W$
as above a $\Z$-bordism for the knot $K$. The following
twisted signatures are associated to it: Define
$\sigma_\omega(W)$ to be the signature of the intersection
form $ \lambda_\omega(W)$ on $H_2(W;
\C_\omega)$, where $\C_\omega$ is the module  over
$\Z[\pi_1W]=\Z[t,t^{-1}]$ obtained by letting $t$ act on $\C$
via multiplication by $ \omega $.  The isomorphism
$\Omega_4(S^1)\cong\Z$, given by the (untwisted!) signature,
implies by the additivity and bordism invariance of all
twisted signatures that the signatures of
$\lambda_\omega(W)$ are in fact independent of the
$\Z$-bordism $W$ for $K$.

The following result was proven in \cite[Section 5]{COT1}
and will be used in order to avoid having to recall the
definition of the $L^2$-signature. Note however that the
(reduced) $L^2$-signature of a knot is by definition the
(reduced) $L^2$-signature of a $\Z$-bordism $W$. By
\cite[Lemma 5.9.4]{COT1} this is independent of the choice of
$W$. Note also that since we have required that the ordinary
signature of $W$ is zero, the reduced and unreduced
$L^2$-signatures are identical.

\begin{lem}\label{lem:5.4}
Let $W$ be a $\Z$-bordism for a knot $K$. Then the
$L^2$-signature of $K$ is the integral  of the twisted
signatures $\sigma_\omega(W)$, integrated over the circle
normalized to length one.
\end{lem}

In order to connect to the setting of \cite[Section 5]{COT1},
note that the universal coefficient spectral sequence implies
that for any $\C[t,t^{-1}]$-module $\mathcal M$ there is an isomorphism
\begin{equation*}
H_2(W;M) \cong H_2(W;\C[t,t^{-1}])
\otimes_{\C[t,t^{-1}]} \mathcal M.
\end{equation*}
First apply this isomorphism to $\mathcal M= \UU\Z$, the algebra of
unbounded operators affiliated to the von Neumann algebra
$\NN\Z$. It follows that the intersection form used in
\cite[Def.5.8]{COT1} to define the $L^2$-signature of $W$ is
nothing else but $\lambda_W\otimes\UU\Z$, where $\lambda_W$ is
Wall's $\Z[t,t^{-1}]$-valued intersection form on
$H_2(W;\Z[t,t^{-1}])$. This implies by
\cite[Lemma 5.6]{COT1} that the $L^2$-signature of $W$ is the
$L^2$-signature of the form $\lambda_W\otimes_{\Z} \C$ on the free
$\C[t,t^{-1}]$-module $H_2(W;\C[t,t^{-1}])/\Tors$.

Finally, apply the isomorphism above to $\mathcal M=\C_\omega$ and
observe that the intersection form $ \lambda_\omega$ is
obtained from the $\C[t,t^{-1}]$-valued intersection form
$\lambda_W\otimes\C$ by substituting $ \omega $ for $t$.
Hence the discussion of \cite[Def.5.3]{COT1} and in
particular Lemma 5.4 apply to prove Lemma~\ref{lem:5.4}.

We return to the proof of Proposition~\ref{LT} above. It now
suffices to prove the following lemma.

\begin{lem} \label{lem:F=W}
For a given Seifert surface $F$ for $K$, there exists a
$\Z$-bordism $W_F$ such that the signatures of
$\lambda_\omega(F)$ and $\lambda_\omega(W_F)$ agree for all
$\omega \in S^1$.
\end{lem}

\begin{proof}
First observe that the equality is true for the
untwisted case $\omega=1$ since both signatures are zero.
This holds for any $\Z$-bordism $W$. All other signatures
can be calculated after inverting the element $(1-t)$ since
for $\omega\neq 1$, the twisted homology is given as follows:
$$
H_2(W;\C_\omega) \cong H_2(W;\C[t,t^{-1}])
\otimes_{\C[t,t^{-1}]} \C_\omega \cong H_2(W; \Lambda)
\otimes_{ \Lambda} \C_\omega
$$
where $ \Lambda: = \C[t,t^{-1}, (1-t)^{-1}]$. The second
isomorphism uses the fact that $ \Lambda$ is a flat module
over the ring
$\C[t,t^{-1}]$. Note that we are using $ \omega \neq 1$ to
make $\C$ a module over $ \Lambda$ (again denoted by
$\C_\omega$).

We finish our proof by showing that for a given Seifert
surface $F$ for $K$, there is a certain choice of $W=W_F$
such that the intersection form on
$H_2(W_F; \Lambda)$ is represented by the matrix
$$
\lambda(F):= (1- t^{-1} ) \cdot S + (1- t)\cdot S^T
$$
where $S$ is again the Seifert matrix for a basis of embedded
curves $a_1,\dots,a_{2g}$ on $F$.
The first part of our computation follows \cite{Ko}
closely.

Let $V_F$ be the complement of a neighborhood of the Seifert surface $F$
pushed into the 4-ball $D^4$.  We will later modify $V_F$ to get the
desired 4-manifold $W_F$.
$\pi_1(V_F) = \langle t \rangle \cong\Z$.  We construct the universal
cover of $V_F$ in the usual manner.  Cut along the trace of pushing the
surface $F$ into the 4-ball.  The cut manifold is homeomorphic to
$D^4$.  This determines two embeddings of
$F\times I$ in the boundary of $D^4$ and we label the positive side
$F_+$ and the negative side $F_-$. Construct the universal cover of $V_F$
from countably many copies of $D^4$ labeled $t^kD^4, k \in \Z$ by
identifying
$F_+ \subset t^kD^4$ with $F_- \subset t^{k-1}D^4$ in the obvious
manner.  Each 4-ball is a fundamental domain of this cover.

The Mayer-Vietoris sequence of this decomposition gives an isomorphism of
$\Z[\Z]$-modules
$$
H_2(V_F; \Z[\Z]) \xrightarrow{\partial} H_1(F; \Z) \otimes_{\Z}
\Z[\Z].
$$
Thus, $H_2(V_F; \Z[\Z])$ is a free $\Z[\Z]$ module generated by $H_1(F)$.
Furthermore, a basis of $H_2(V_F; \Z[\Z]) \cong \pi_2(V_F)$ is obtained from a basis
of curves of $H_1(F;\ Z)$,  $a_i \subset F$, by choosing immersed 2-disks 
$$
(\Sigma_i^+,\partial\Sigma_i^+) \subset (tD^4,F_+) \text{ and }
(\Sigma_i^-,\partial\Sigma_i^-) \subset (t^0D^4,F_-)
$$
and orienting $\Sigma_i := \Sigma_i^+ \cup \Sigma_i^-$ (and $F$) so that 
$$
\partial\Sigma_i^+=a_i^\uparrow \text{ and } \partial\Sigma_i^-= -a_i^\downarrow.
$$
By the Mayer-Vietoris argument given above and the usual geometric interpretation of
the connecting homomorphism, the spheres $\Sigma_1,
\ldots, \Sigma_{2g}$ together give a basis for the free $\Z[\Z]$ module $H_2(V_F;
\Z[\Z])$.  We now compute the equivariant intersection form $\lambda_{V_F}$ on
$H_2(V_F; \Z[\Z]) \cong \pi_2(V_F)$ given by
$$
\lambda_{V_F} (\Sigma_i, \Sigma_j) = \sum_k (\Sigma_i \cdot
t^k \Sigma_j) t^k \quad\in \Z[\Z],
$$
where $\cdot$ is the usual (integer valued) intersection number of the $2$-spheres
$\Sigma_i$ in the universal cover of $V_F$.

By general position, we may assume that $\Sigma_i$ and $\Sigma_j$ do not
intersect each other along $F\times I$ (since they intersect $F\times I$ only in their
boundary circles). 
Moreover, $\Sigma_i \cdot t^k \Sigma_j = 0$ for $k \neq -1,0,1$. Recalling
that $S_{i,j}= lk(a_i, a_j^\uparrow)$ is the $(i,j)$ entry of the Seifert matrix
$S$ associated to the Seifert surface $F$, we now compute equivariant intersection
numbers as follows.
\begin{eqnarray*}
\Sigma_i \cdot \Sigma_j &=& \Sigma_i^- \cdot \Sigma_j^- + \Sigma_i^+
\cdot \Sigma_j^+  = lk(-a_i^\downarrow, -a_j^\downarrow) +
lk(a_i^\uparrow,a_j^\uparrow) \\
&=& lk(a_i, a_j^\downarrow) + lk(a_i,a_j^\uparrow) = S_{j,i} + S_{i,j}.
\end{eqnarray*}
Note that by symmetry, this is independent of how the curves $a_i$
respectively $a_j$ are pushed into general position. Moreover, one has
\begin{eqnarray*}
\Sigma_i \cdot t\Sigma_j &=& \Sigma_i^+ \cdot t\Sigma_j^- 
= lk(a_i^\uparrow, -a_j^\downarrow) =
-lk(a_j,a_i^\uparrow) = -S_{j,i},\\
\Sigma_i \cdot t^{-1}\Sigma_j &=& \Sigma_i^- \cdot
t^{-1}\Sigma_j^+ =  lk(-a_i^\downarrow, a_j^\uparrow)  =
-lk(a_i, a_j^\uparrow) = -S_{i,j}.
\end{eqnarray*}
Thus the intersection form on $H_2(V_F; \Z[\Z])$ is given by the matrix
$$
\lambda_{V_F}(\Sigma_i,\Sigma_j) =
(1-t^{-1}) S_{i,j} + (1-t) S^T_{i,j} =\lambda(F)_{i,j} 
$$
as claimed.

Even though it might look now like we are done, the
$4$-manifold $V_F$ isn't quite what we want since it's
boundary is {\em not} $M_K$. By construction, $\partial V_F$
is obtained from the knot complement by attaching $F \times
S^1$ (instead of $D^2\times  S^1$ to get $M_K$.)  We turn $V_F$ into a
$4$-manifold $W_F$ with the correct boundary by attaching {\em round
handles} $(D^2
\times D^1) \times S^1$ to a disjointly embedded half basis
of curves on $F$. In other words, we start with $g$ surgeries
turning $F$ into a 2-disk and then cross with the circle. Clearly the
boundary of $W_F$ is $M_K$.  What happened to the
equivariant intersection form?

Here we make use of the ring $\Lambda$ to simplify life
tremendously.  Lemma~\ref{lem:Lambda} below and an easy Mayer-Vietoris
sequence argument imply that the inclusion
induces an isomorphism $H_2(V_F; \Lambda) \cong H_2(W_F;\Lambda)$, and
thus concludes the proof of Proposition~\ref{LT}.
\end{proof}

\begin{rem}
The above proof inverts $1-t$ to save a lot of work.  The more
geometric reader will
find the intersection pairing
for $W_F$ {\em before localization} computed in
\cite{Ko}.
\end{rem}

One should compare that work with the simplicity of the
following

\begin{lem}\label{lem:Lambda} Let $X$ be a finite CW-complex
and $t$ be a generator of $\pi_1S^1$. If $\Lambda: =
\Z[t,t^{-1}, (1-t)^{-1}]$, then
$$
H_*(X \times  S^1; \Lambda) = 0
$$
if the twisting is given by the projection $p_2:X \times
S^1\to S^1$.
\end{lem}

\begin{proof}
The differential in the cellular chain complex of the universal
cover of $S^1=e^0\cup e^1$ is given by multiplication with
$(t-1)$. This is clearly an isomorphism after tensoring with
$\Lambda$.  Hence, this chain complex is acyclic, and thus contractible.

The cellular chain complex of the $\Z$-cover of $X\times S^1$ with
$\Lambda$ coefficients
is the tensor product of
the cellular chains of $X$ with the contractible
$\Lambda$-module chain complex for $S^1$.  A tensor product with a
contractible chain complex is again contractible, and therefor acyclic.
\end{proof}

\begin{rem} It is clear that the manifold $V_F$ constructed in
Lemma~\ref{lem:F=W} has a spin structure as a subset of the
4-ball. The reader is invited to check that the $\Z$-bordism
$W_F$ is a spin manifold if and only if the original knot $K$
has trivial Arf invariant. For example, if the Arf invariant
is trivial then one may choose a half basis of curves
$a_1,\dots,a_g$ on $F$ with even self-linking number. This
implies that the round surgeries leading from $V_F$ to $W_F$
are spin structure preserving.
\end{rem}

We now return to prove the technical Proposition
\ref{independence}.

\begin{proof}[Proof of 2.6] Let $J_m$, $m\ge1$, be a knot
with Alexander polynomial $\Delta(t)=2mt-(4m-1)+2mt^{-1}$.
Such a knot necessarily has zero Arf invariant (see
Theorem~10.4 of \cite{Ka}). Let $M$ be the zero surgery on
$J_m$ and $W$ be a compact, spin $4$-manifold with $\p W=M$,
$\pi_1(W)\cong\Z$ and $j_*:H_1(M;\Z)\lra H_1(W,\Z)$ an
isomorphism. By considering the long exact homology sequence
for the pair $(W,M)$ with $\Z[t,t^{-1}]$ coefficients (using
the canonical epimorphisms $\phi_i$, $\psi_i$ to $\Z$ to
define the coefficient systems) one sees that the order of
the Alexander module $H_1(M;\Q[t,t^{-1}])$ is the determinant
of a matrix $h$ representing the intersection form on
$H_2(W;\Q[t,t^{-1}])$. The (reduced) signature function
$\s(h(z))-\s_0(h):S^1\lra\Z$ is a locally-constant function
which is $0$ at $z=1$ and changes value only (possibly) at
the two zeros of $\det h(z)=\pm\Delta(z)$ which are $e^{\pm
i\tha_m}$, where $0\le\tha_m\le\f\pi2$ and
$\cos\tha_m=(4m-1)/4m$. Since $\det h(z)$ changes sign at
$e^{i\tha_m}$, $\s(h(z))$ {\it must} change value there. Thus
$\s(h(-1))=a_m$, for some non-zero integer $a_m$ (actually
$\pm2$), and $\rho(M_m,\phi_m)=a_m(\pi-\tha_m)/\pi$. To prove
2.6, it will suffice to prove that there exists a infinite
collection of integers $m$ such that $\{\tha_m\}$ is linearly
independent over the integers.

Choose primes $p_i$, $5\le p_1<p_2<\dots<p_j<\dots$, each
congruent to $1$ modulo~$4$, and set $m_j=\f18(p_j-1)^2$. We
claim that if $\Th_j=\cos^{-1}\Bigl(\f{4m_j-1}{4m_j}\Bigr)$,
then $\{\Th_j\}$ is linearly independent over the
rationals.

\begin{lem}\label{galois} Suppose $5\le p_1<p_2<\dots<p_n$
are primes. Let $\xi_j=i\sqrt{p_j(p_j-2)}$ where
$i=\sqrt{-1}$. Then $[\Q(\xi_1,\dots,\xi_n):\Q]=2^n$ and
the Galois group over $\Q$ is $(\Z_2)^n$ generated by
automorphisms $\phi_j$ where
$\phi_j(\xi_\ell)=\begin{cases}-\xi_\ell&
\text{if }\ \ell=j\\
\xi_\ell& \text{if }\ \ell\neq j
\end{cases}$.
\end{lem}

\begin{proof}[Proof of Lemma \ref{galois}] One verifies
easily that the Lemma is true for $n=1$. Now assume
$[\Q(\xi_1,\dots,\xi_{k-1}):\Q]=2^{k-1}$ for {\it any}
increasing sequence of primes $5<p_1<\dots<p_{k-1}$. We shall
show $[\Q(\xi_1,\dots,\xi_k):\Q]=2^k$. It suffices to show
$\xi_k$ is not in $\Q(\xi_1,\dots,\xi_{k-1})$. If it were
then there are $\al$, $\beta\in\Q(\xi_1,\dots,\xi_{k-2})$ with
$\xi_k=\al+\beta\xi_{k-1}$. Squaring each side yields
$$
\xi^2_k = (\al^2 + \beta^2\xi^2_{k-1}) + 2\al\beta\xi_{k-1}.
$$
Since $\xi^2_k$, $\xi^2_{k-1}$, $\al$ and $\beta$
lie in $\Q(\xi_1,\dots,\xi_{k-2})$ and since, by the inductive hypothesis,
$\xi_{k-1}\notin\Q(\xi_1,\dots,\xi_{k-2})$, we have
$\al\beta=0$. If $\beta=0$ then
$\xi_k\in\Q(\xi_1,\dots,\xi_{k-2})$. This will contradict the
induction hypothesis for the sequence
$5<p_1<\dots<p_{k-2}<p_k$. Suppose $\al=0$. Since
$\beta=\g+\de\xi_{k-2}$ where $\g$,
$\de\in\Q(\xi_1,\dots,\xi_{k-3})$, by squaring each side of
$\xi_k=\beta\xi_{k-1}$ and applying induction we again
conclude $\g\de=0$. As above if $\de=0$ then
$\beta\in\Q(\xi_1,\dots,\xi_{k-3})$ so
$\xi_k\in\Q(\xi_1,\dots,\xi_{k-3},\xi_{k-1})$ contradicting
the induction hypothesis. Thus $\g=0$ and
$\xi_k=\de\xi_{k-2}\xi_{k-1}$. We continue in this fashion
until we conclude $\xi_k=\Th\xi_1\xi_2\dots\xi_{k-1}$ where
$\Th\in\Q$, implying
$$
-p_k(p_k - 2) = (r^2/s^2)\xi^2_1\dots\xi^2_{k-1}.
$$
Multiplying by $s^2$ yields an equation over $\Z$ which has
no solution since $p_k$ divides the left hand side an odd
number of times and the right hand side an even number of
times. This contradiction finishes the proof that
$[\Q(\xi_1,\dots,\xi_k):\Q]=2^k$.

Consider the subgroup of the Galois group generated by the
$\phi_j$. The subset
$\{\phi_{i_1}\phi_{i_2}\dots\phi_{i_m}\mid1\le\phi_{i_1}<\dots<\phi_{i_m}\le
n\}$ is easily seen to contain $2^n$ distinct elements of
order 2 and the result follows.
\end{proof}

Now suppose $\sum^n_{j=1}c_j\Th_j=0$ where $c_j\in\Z$.
Multiplying by $i$ and exponentiating yields
$$
\prod^n_{j=1}\Biggl(\f{4m_j-1}{4m_j} +
\f{\xi_j}{4m_j}\Biggr)^{c_j} = 1.
$$
Evaluating $\phi_j$ on each side
yields that
$\Bigl(\f{4m_j-1}{4m_j}+\f{\xi_j}{4m_j}\Bigr)^{c_j}$ is
real and hence (being of norm~$1$) is $\pm1$. Thus
$\f{4m_j-1}{4m_j}+\f{\xi_j}{4m_j}$ is a primitive
$r^{\text{th}}_j$ root of unity for some $r_j\mid c_j$. Since
$[\Q(\xi_j),\Q]=2$, $r_j$ equals $3$, $4$ or $6$. But, given
the bound on $m_j$, one checks by inspection that this is not
the case. This is a contradiction.
\end{proof}

\end{document}